\documentclass{article}
\usepackage{amsmath}
\usepackage{amssymb}

\def\be{\begin{equation}}       \def\ee{\end{equation}}
\def\bd{\begin{displaymath}}    \def\ed{\end{displaymath}}
\def\beq{\begin{eqnarray}}      \def\eeq{\end{eqnarray}}
\def\bseq{\begin{eqnarray*}}    \def\eseq{\end{eqnarray*}}
\def\ba{\begin{array}}          \def\ea{\end{array}}
\def\ben{\begin{enumerate}}     \def\een{\end{enumerate}}

\newtheorem{Prop}{Proposition}
\def\bpr{\begin{Prop}}           \def\epr{\end{Prop}}

                     \def\cl{{\cal L}}
                     
                    \def\ct{{\cal T}}

\def\ta{\tilde{A}}
\def\tb{\tilde{B}}
\def\tc{\tilde{C}}
\def\td{\tilde{D}}
\def\tf{\tilde{F}}

\def\lra{\longrightarrow}

\def\cop{\Delta}

\def\cnt{\varepsilon}
\def\ot{\otimes}

\def\GLqtwo{{GL}_{q}(2)}
\def\GLqthree{{GL}_{q}(3)}

\def\GLpqtwo{{GL}_{p,q}(2)}
\def\GLqtwoC{{GL}_{q}(2,C)}

\def\Grs{{G}_{r,s}}
\def\GLone{{GL}({1})}
\def\GLtwo{{GL}({2})}
\def\GLthree{{GL}({3})}
\def\GLn{{GL}({n})}
\def\rinv{{r}^{-1}}
\def\sinv{{s}^{-1}}
\def\Casmir{{\bf D}}
\def\Casmirinv{{\Casmir}^{-1}}
\def\ident{{\bf 1}}

\begin{document}

\begin{center}
{\large\bf 
DUALITY FOR THE $\Grs$ QUANTUM GROUP\footnote{Kyoto University preprint
RIMS - 1260 (1999).}
}

\bigskip
\bigskip
\bigskip
\bigskip
\bigskip

{\large\bf Deepak Parashar and Roger J. McDermott
}

\bigskip
School of Computer and Mathematical Sciences\\
The Robert Gordon University\\
St. Andrew Street\\
Aberdeen AB25 1HG\\
United Kingdom\\

{\sl 
email: deeps@scms.rgu.ac.uk, rm@scms.rgu.ac.uk
}

\bigskip
\bigskip
\bigskip
\bigskip
\bigskip
{\bf Abstract}
\end{center}
\medskip

The two parameter quantum group $G_{r,s}$ is generated by five elements,
four of which form a Hopf subalgebra isomorphic to $GL_{q}(2)$, while the
fifth generator relates $G_{r,s}$ to $GL_{p,q}(2)$. We construct
explicitly the dual algebra of $G_{r,s}$ and show that it is isomorphic to
the single parameter deformation of $gl(2)\oplus gl(1)$, with the second
parameter appearing in the costructure. We also formulate a
differential calculus on $\Grs$ which provides a realisation of the
calculus on $\GLpqtwo$.\par


\newpage
\section{Introduction}
The q-deformation of the universal enveloping algebra is commonly referred
to in recent literature as quantum universal enveloping algebras or simply
quantum algebras. There are some approaches [1-5] to quantum groups
in which the objects may be called {\em quantum matrix groups} which are
Hopf algebras in duality to the quantum algebras. In particular this means
the existence of a doubly
nondegenerate bilinear form between the two Hopf algebras. It was Sudbery
[6] who first gave a formalism for such a duality motivated by the fact
that, at the classical level, an element of the Lie algebra corresponding
to a Lie group is a tangent vector at the identity of the Lie group. The
q-analogue of tangent vector at the identity would then be obtained by
first differentiating the elements of the given Hopf algebra H and then
taking the counit operation. The elements thus obtained would belong to
the dual Hopf algebra H*. Following this simple procedure, Sudbery
obtained $U_{q}(sl(2))\ot U(u(1))$ as the algebra of tangent vectors at
the identity of $\GLqtwo$. This technique has been successfully applied to
the case of multiparametric $\GLpqtwo$ [7] and quantum $\GLn$ [8]
groups. The quantised universal enveloping algebra (i.e. the dual)
corresponding to a quantum group can also be obtained using the $R$-matrix
formalism [1].

\par

Recently, there has been a considerable development in the study of
multiparametric quantum groups, both from algebraic and differential
geometric aspects. We wish to focus our attention on a new quantum group
[9] $\Grs$, depending on two deformation
parameters and five generators. The first four generators of this Hopf
algebra form a Hopf subalgebra, which coincides exactly with the single
parameter dependent $\GLqtwo$ quantum group when $q = \rinv$.  However,
the two parameter dependent $\GLpqtwo$ can also be realised through the
generators of this $\Grs$ Hopf algebra, provided the sets of deformation
parameters $p$, $q$ and $r$, $s$ are related to each other in a
particular fashion. This new algebra can, therefore, be used to realise
both $\GLqtwo$ and $\GLpqtwo$ quantum groups. Alternatively, this $\Grs$
structure can be considered as a two parameter quantisation of the
classical $\GLtwo\times \GLone$ group.  The first four generators of
$\Grs$, i.e. $a$, $b$, $c$, $d$ correspond to $\GLtwo$ group at the
classical level and the remaining generator $f$ is related to $\GLone$
group. Given the rich structure of this simple quantum group it is
surprising that no further studies have so far been done. A starting point
in the further investigation would be the explicit description of the dual
algebra to $\Grs$ and to look at its differential structure.\par

This is the problem that we address in this paper. We obtain the dual
algebra to the two-parameter matrix quantum group $\Grs$.
We show that the Hopf algebra dual to $\Grs$ may be realised using the
method described by Sudbery [6] and exhibit its Hopf Structure. The dual
algebra is also explicitly constructed using the $R$-matrix approach [1].
We then employ Jurco's constructive procedure [10] for a description of
the differential calculus on $\Grs$.

In Section II, we describe the essential features of $\Grs$. The dual
algebra for $\Grs$ is then constructed in Section III within the framework
of Sudbery's approach wherein we also detail the commutation relations
and Hopf structure of the generators of the dual. Section IV deals
with the $R$-matrix formulation of $\Grs$ while Section V is devoted to
the construction of a bicovariant differential calculus on $\Grs$.
Concluding remarks are made in Section VI.

\section{The quantum group $\Grs$}

The Hopf algebra $\Grs$ is generated by elements $a$, $b$, $c$, $d$, and
$f$ satisfying the relations
\[
\begin{array}{ll}
ab=\rinv ba,&db=rbd\\
ac=\rinv ca,&dc=rcd\\
bc=cb,&[a,d]=(\rinv-r)bc
\end{array}
\]

and
\[
\begin{array}{ll}
af=fa,&cf=sfc\\
bf=\sinv fb,&df=fd
\end{array}
\]

where $r$ and $s$ are two deformation parameters with arbitrary nonzero
complex numbers.  Elements $a$, $b$, $c$, $d$ satisfying the first set of
commutation relations form a subalgebra which coincides exactly with
$\GLqtwo$ when $q = \rinv$.  The generators are arranged in a $3 \times
3$ matrix as
\[
T=\begin{pmatrix}a&b&0\\c&d&0\\0&0&f\end{pmatrix}
\]

The Hopf structure is given as
\[
\begin{array}{l}
\cop\begin{pmatrix}a&b&0\\c&d&0\\0&0&f\end{pmatrix} = 
\begin{pmatrix}
a\ot a+b\ot c & a\ot b+b\ot d & 0\\
c\ot a+d\ot c & c\ot b+d\ot d & 0\\
0 & 0 & f\ot f
\end{pmatrix}\\
\\
\cnt\begin{pmatrix}a&b&0\\c&d&0\\0&0&f\end{pmatrix} = 
\begin{pmatrix}1&0&0\\0&1&0\\0&0&1\end{pmatrix}
\end{array}
\]

The Casimir operator is defined as $\Casmir = ad-\rinv bc$. The inverse is
assumed to exist and satisfies $\cop(\Casmirinv)=\Casmirinv\ot\Casmirinv$, $\cnt(\Casmirinv)=1$, $S(\Casmirinv)=\Casmir$, which enables one to determine the antipode matrix $S(T)$, as
\[
S\begin{pmatrix}a&b&0\\c&d&0\\0&0&f\end{pmatrix}
=\Casmirinv
\begin{pmatrix}d&-rc&0\\-\rinv c&a&0\\0&0&\Casmir f\end{pmatrix}
\]

The elements $a$, $b$, $c$ and $d$ of $\Grs$ evidently form a Hopf
subalgebra, which coincides with $\GLqtwo$ quantum group.  In order to
investigate the connection of $\Grs$ with $\GLpqtwo$, a simple realisation
of $\GLpqtwo$ generators through the elements of $\Grs$ has been proposed as 
\[
\begin{pmatrix}a'&b'\\c'&d'\end{pmatrix}
=f^{N}
\begin{pmatrix}a&b\\c&d\end{pmatrix}
\]

where $N$ is any fixed nonzero integer.  For this to give a realisation of $\GLpqtwo$, the deformation parameters $(p,q)$ and $(r,s)$ are related via
\[p = \rinv s^{N} \quad \text{and} \quad q = \rinv s^{-N}\]

For $N\neq 0$ and any given value of $(p,q)$, one can find out the corresponding values of the deformation parameters $(r,s)$.  The full Hopf algebra structure related to $\GLpqtwo$ can in fact be reproduced through the above realisation.  The mapping from $(r,s)$-plane to $(p,q)$-plane depends on the choice of $N$ $(\neq 0)$.  Taking different values of $N$, a single point on the $(r,s)$-plane can be mapped over infinite number of discrete points which satisfy
\[
p_{N}q_{N}=r^{-2},\quad
\frac{p_{N+1}}{p_{N}}=s,\quad
\frac{q_{N+1}}{q_{N}}=\sinv,
\]

where
\[ p_{N} = \rinv s^{N}\quad    \text{and} \quad q_{N} = \rinv s^{-N} \]

Here $(p_{N}, q_{N})$ denotes a point on the $(p,q)$-plane corresponding
to the $N$-th mapping.  One finds that from the representation of $\Grs$
algebra for a particular value of deformation parameters $(r,s)$, one can
build up the representations of $\GLpqtwo$ quantum group for infinitely
many discrete values of $(p,q)$ parameters.  The choice of $(9 \times 9)$
R-matrix is given by
\[
R=\sum_{i}{e_{ii}\ot e_{ii}} + \sum_{i\neq j}{f_{ii}\cdot e_{ii}\ot
e_{jj}} + (r-\rinv)\sum_{i<j}{e_{ij}\ot e_{ji}}
\]

where $f_{12} = f_{23} = 1$, $f_{13} = s$ and $f_{ij} = f^{-1}_{ji}$ and
$i, j \in [1,3]$.  This satisfies the quantum Yang-Baxter equation which
when recast in the form of the RTT-relations yields the commutators of
the elements $\Grs$.  In fact, the most
general Hopf algebra generated by this R-matrix is $\GLqthree$ which can
easily reproduce the $\Grs$ structure. Therefore, $\Grs$ might be
interpreted as a quotient of multi-parameter deformed $\GLthree$
group.\par

\section{The Dual algebra for $\Grs$}

Two bialgebras $U$ and $A$ are in duality if there exists a doubly
non-degenerate bilinear form
\begin{equation*}
\langle,\rangle:U\ot A\rightarrow C;\quad
\langle,\rangle:(u,a)\rightarrow \langle u,a\rangle;\quad
\forall u\in U, a\in A
\end{equation*}

such that for $u,v\in U$ and $a,b\in A$, have
\begin{equation*}
\begin{array}{ll}
\langle u,ab\rangle = \langle \cop_{U}(u), a\ot b\rangle\\
\langle uv,a\rangle = \langle u\ot v, \cop_{A}(a)
\end{array}
\end{equation*}

\begin{equation*}
\begin{array}{ll}
\langle \ident_{U},a\rangle=\cnt_{A}(a)\\
\langle u,\ident_{A}\rangle=\cnt_{U}(u)
\end{array}
\end{equation*}

For the two bialgebras to be in duality as Hopf algebras, $U$ and $A$
further satisfy
\begin{equation*}
\langle S_{U}(u),a \rangle=\langle u,S_{A}(a)\rangle
\end{equation*}

It is enough to define the pairing between the generating elements of
the two algebras. Pairing for any other elements of $U$ and $A$ follows
from the above relations and the bilinear form inherited by the
tensor product. For example, for $\cop (u)=\sum_{i}{{u'}_{i}\ot
{u''}_{i}}$, we have
\[
\begin{align*}
\langle u,ab\rangle = \langle \cop_{U}(u),a\ot b\rangle &
  =\sum_{i}{\langle {u'}_{i}\ot {u''}_{i}, a\ot b\rangle}\\
\quad&
  =\sum_{i}{\langle {u'}_{i},a\rangle\langle {u''}_{i}, b\rangle}
\end{align*}
\]

Sudbery's idea for duality was motivated by the fact that, at the
classical level, "an element of Lie Algebra corresponding to the Lie
Group is a tangent vector of the identity of the Lie Group". Let
$H$ be a given Hopf algebra generated by non-commuting elements $a$,
$b$, $c$, $d$. Then $q$-analogue of tangent vector of identity is
obtained by differentiating elements of $H$ (polynomials in $a$, $b$,
$c$, $d$) and then putting  
$\left( \begin{smallmatrix}a&b\\c&d\end{smallmatrix} \right) =
\left( \begin{smallmatrix}1&0\\0&1\end{smallmatrix} \right)$
later on (i.e. taking the co-unit operation). The elements thus obtained
belong to the dual $H^{*}$. For the sake of convenience differentiation is
taken from the right. He obtained 
$U_{q}(sl(2))\ot U(u(1))$ as the algebra of tangent vectors at the
identity of $\GLqtwoC$.\par

As a Hopf algebra, $\Grs$ is generated by elements $a$, $b$, $c$, $d$, $f$
and a basis is given by all monomials of the form
\[ g = g_{kltmn} = a^{k}d^{l}f^{t}b^{m}c^{n} \]

where  
$k,l,t,m,n\in Z_{+}$, and $\delta_{00000}$ is the unit of the algebra 
$\ident_{A}$. We use the so-called normal ordering i.e. first put the
diagonal elements from the T-matrix then use the lexicographic order for
the others. Let $U_{r,s}$ be the algebra generated by tangent vectors at
the identity of $\Grs$. Then $U_{r,s}$ is dually paired with $\Grs$. The
pairing is defined through the tangent vectors as follows
\[
\langle Y,g\rangle = \frac{\partial g}{\partial y}|_{
  \left(      
  \begin{smallmatrix}a&b&0\\c&d&0\\0&0&f\end{smallmatrix}\right) =
  \left( \begin{smallmatrix}1&0&0\\0&1&0\\0&0&1\end{smallmatrix}\right)
}\quad;\quad
y\in\left\{ 
a,b,c,d,f
\right\}
\]

or more consistently
\[
\langle Y,g\rangle = \cnt(\frac{\partial g}{\partial y})
\quad;\quad
y\in\left\{ 
a,b,c,d,f
\right\}
\]
where $Y\in \left\{ A,B,C,D,F\right\}$ are the generating elements of the
dual algebra (which has unit $\ident_{U}$). Such a generalised approach in
which initial pairings are postulated to be equal to the classical
undeformed results was given by Dobrev in [7]. Writing explicitly the 5
possible cases, we have
\begin{eqnarray*}
\langle A,g\rangle &=& \cnt(\frac{\partial g}{\partial a}) 
  = k\delta_{m0}\delta_{n0}\\
\langle B,g\rangle &=& \cnt(\frac{\partial g}{\partial b})
   = \delta_{m1}\delta_{n0}\\
\langle C,g\rangle &=& \cnt(\frac{\partial g}{\partial c})
   = \delta_{m0}\delta_{n1}\\
\langle D,g\rangle &=& \cnt(\frac{\partial g}{\partial d})
   = l\delta_{m0}\delta_{n0}\\
\langle F,g\rangle &=& \cnt(\frac{\partial g}{\partial f})
   = t\delta_{m0}\delta_{n0}
\end{eqnarray*}

where differentiation is from the right as this is most suitable
for differentiation in this basis. Variables $b$, $c$ commute and their
differentiation is classical. For example, $\langle B,g\rangle$ is
obtained by first calculating 
\[
(a^{k}d^{l}f^{t}b^{m}c^{n})\frac{\partial}{\partial b}=
(a^{k}d^{l}f^{t}c^{n}b^{m})\frac{\partial}{\partial b}=
a^{k}d^{l}f^{t}c^{n}mb^{m-1}
\]
	
(because of differentiation from right), and then applying the counit
operation
\[
\cnt(a^{k}d^{l}f^{t}c^{n}mb^{m-1}) = 
{\cnt(a)}^{k} {\cnt(d)}^{l} {\cnt(f)}^{t} {\cnt(c)}^{n} m {\cnt(b)}^{m-1} =
\delta_{n0}m\delta_{m-1,0}=\delta_{n0}\delta_{m1}
\]

Other pairings follow in a similar way. As a consequence of the above
pairings, the following relations hold
\begin{eqnarray*}
\langle A, \left( 
\begin{smallmatrix}
a&b&0\\c&d&0\\0&0&f
\end{smallmatrix}
\right) \rangle
= \left( \begin{smallmatrix}
1&0&0\\0&0&0\\0&0&0
\end{smallmatrix}\right)
\end{eqnarray*}
\begin{eqnarray*}
\langle D, \left( 
\begin{smallmatrix}
a&b&0\\c&d&0\\0&0&f
\end{smallmatrix}
\right) \rangle
= \left( \begin{smallmatrix}
0&0&0\\0&1&0\\0&0&0
\end{smallmatrix}\right)
\end{eqnarray*}
\begin{eqnarray*}
\langle F, \left( 
\begin{smallmatrix}
a&b&0\\c&d&0\\0&0&f
\end{smallmatrix}
\right) \rangle
= \left( \begin{smallmatrix}
0&0&0\\0&0&0\\0&0&1
\end{smallmatrix}\right)
\end{eqnarray*}
\begin{eqnarray*}
\langle B, \left( 
\begin{smallmatrix}
a&b&0\\c&d&0\\0&0&f
\end{smallmatrix}
\right) \rangle
= \left( \begin{smallmatrix}
0&1&0\\0&0&0\\0&0&0
\end{smallmatrix}\right)
\end{eqnarray*}
\begin{eqnarray*}
\langle C, \left( 
\begin{smallmatrix}
a&b&0\\c&d&0\\0&0&f
\end{smallmatrix}
\right) \rangle
= \left( \begin{smallmatrix}
0&0&0\\1&0&0\\0&0&0
\end{smallmatrix}\right)
\end{eqnarray*}
\[
\langle Y, \ident_{A}\rangle = 0\quad Y=A,B,C,D,F
\]
\[
\langle\ident_{u},a^{k}d^{l}f^{t}b^{m}c^{n}\rangle = \delta_{m0}\delta_{n0}
\]
\[
\left[ \because LHS = \cnt_{A}(a^{k}d^{l}f^{t}b^{m}c^{n}) =
{(\cnt_{A}(a))}^{k} {(\cnt_{A}(d))}^{l} {(\cnt_{A}(f))}^{t}
{(\cnt_{A}(b))}^{m} {(\cnt_{A}(c))}^{n} =\delta_{m0}\delta_{n0}
\right]
\]
%\newpage
\subsection*{Commutation Relations}
The action of the monomials in $U_{r,s}$ on $g =
a^{k}d^{l}f^{t}b^{m}c^{n}$ is given by the following:

\begin{eqnarray*}
\langle BC, g\rangle &=&
\delta_{m0}\delta_{n0}\sum_{j=0}^{k-1}{r^{2(j-l)}}
+ r^{-1}\delta_{m1}\delta_{n1}\\
\langle CB, g\rangle &=&
\delta_{m0}\delta_{n0}\sum_{j=0}^{k-1}{r^{-2j}}
+ r\delta_{m1}\delta_{n1}\\
\langle AB, g\rangle &=& (k+1)\delta_{m1}\delta_{n0} = (k+1)\langle B,
g\rangle\\
\langle BA, g\rangle &=& k\delta_{m1}\delta_{n0} = k\langle B, g\rangle\\
\langle AC, g\rangle &=& k\delta_{m0}\delta_{n1} = k\langle C, g\rangle\\
\langle CA, g\rangle &=& (k+1)\delta_{m0}\delta_{n1} = (k+1)\langle C,
g\rangle\\
\langle DB, g\rangle &=& l\delta_{m1}\delta_{n0} = l\langle B, g\rangle\\
\langle BD, g\rangle &=& (l+1)\delta_{m1}\delta_{n0} = (l+1)\langle B,
g\rangle\\
\langle DC, g\rangle &=& (l+1)\delta_{m0}\delta_{n1} = (l+1)\langle C,
g\rangle\\
\langle CD, g\rangle &=& l\delta_{m0}\delta_{n1} = l\langle C, g\rangle\\
\langle AD, g\rangle &=& \langle DA, g\rangle = kl\delta_{m0}\delta_{n0} =
kl\langle\ident_{U}, g\rangle\\
\langle AF, g\rangle &=& \langle FA, g\rangle = kt\delta_{m0}\delta_{n0} =
kt\langle\ident_{U}, g\rangle\\
\langle DF, g\rangle &=& \langle FD, g\rangle = lt\delta_{m0}\delta_{n0} =
lt\langle\ident_{U}, g\rangle\\
\langle BF, g\rangle &=& (t+1)\delta_{m1}\delta_{n0} = (t+1)\langle B,
g\rangle\\
\langle FB, g\rangle &=& t\delta_{m1}\delta_{n0} = t\langle B, g\rangle\\
\langle CF, g\rangle &=& t\delta_{m0}\delta_{n1} = t\langle C, g\rangle\\
\langle FC, g\rangle &=& (t+1)\delta_{m0}\delta_{n1} = (t+1)\langle C,
g\rangle
\end{eqnarray*}

\newpage
As an illustration, $\langle BC, g\rangle$ is obtained by using the
assumed duality,
\begin{eqnarray*}
&\langle BC,& a^{k}d^{l}f^{t}b^{m}c^{n}\rangle\\
&=& \langle B\ot C, \cop_{A}(a^{k}d^{l}f^{t}b^{m}c^{n})\rangle\\
&=& \langle B\ot C,
(\cop_{A}(a))^{k}(\cop_{A}(d))^{l}(\cop_{A}(f))^{t}(\cop_{A}(b))^{m}
(\cop_{A}(c))^{n}\rangle\\
&=& \langle B\ot C,
(a\ot a+b\ot c)^{k}(c\ot b+d\ot d)^{l}(f\ot f)^{t} (a\ot b+b\ot
d)^{m}(c\ot a+d\ot c)^{n}\rangle\\
\\
&=& \langle B\ot C, (a^{k}\ot a^{k} + \sum_{j=0}^{k-1}{a^{k-1-j}ba^{j}\ot
a^{k-1-j}ca^{j}})(d^{l}\ot d^{l} + \sum_{j=0}^{l-1}{d^{l-1-j}cd^{j}\ot 
d^{l-1-j}bd^{j}})\\
& &(f^{t}\ot f^{t})(\delta_{m0} + (a\ot b+b\ot d)\delta_{m1})(\delta_{n0}
+ (c\ot a+d\ot c)\delta_{n1})\rangle\\
& & \text{(keeping terms involving $b$,$c$ at most of degree 1)}\\
\\
&=&\langle B\ot C,(\sum_{j=0}^{k-1}{a^{k-1-j}ba^{j}d^{l}f^{t}\ot
a^{k-1-j}ca^{j}d^{l}f^{t}})\delta_{m0}\delta_{n0} +
(a^{k}d^{l}f^{t}bd\ot a^{k}d^{l}f^{t}dc)\delta_{m1}\delta_{n1}\\
& & \text{(keeping only contributing terms involving exactly one factor
$'b'$ to the left of $\ot$ and}\\
& & \text{one factor $'c'$ to the right of $\ot$)}
\end{eqnarray*}

Next, use $G_{r,s}$ commutation relations to reorder the elements to
the basis order ($adfbc$). Then, we have
\begin{eqnarray*}
&\langle BC,& a^{k}d^{l}f^{t}b^{m}c^{n}\rangle\\
&=&\langle B\ot C,(\sum_{j=0}^{k-1}{a^{k-1}d^{l}f^{t}br^{j-l}s^{-t}\ot
a^{k-1}d^{l}f^{t}cr^{j-l}s^{t}})\delta_{m0}\delta_{n0}\\
& & + (a^{k}d^{l+1}f^{t}br^{-1}\ot
a^{k}d^{l+1}f^{t}c)\delta_{m1}\delta_{n1}\rangle\\
& & \text{(Moving the elements to the right (left) brings in the powers of
$r$ and $s$)}
\end{eqnarray*}

Now, applying the pairings proves $\langle BC, g\rangle$.  Similarly
for $\langle CB, g\rangle$,we have

\begin{eqnarray*}
&\langle CB,& a^{k}d^{l}f^{t}b^{m}c^{n}\rangle\\
&=& \langle C\ot B, (a^{k}\ot a^{k})(d^{l}\ot d^{l}
+
\sum_{j=0}^{l-1}{d^{l-1-j}cd^{j}\ot d^{l-1-j}bd^{j}})(f^{t}\ot f^{t})\\
& &(\delta_{m0} + (a\ot b)\delta_{m1})(\delta_{n0}+(c\ot 
a)\delta_{n1})\rangle\\
\\
&=&\langle C\ot B, (\sum_{j=0}^{l-1}{a^{k}d^{l-1-j}cd^{j}f^{t}\ot
a^{k}d^{l-1-j}bd^{j}f^{t}})\delta_{m0}\delta_{n0} + (a^{k}d^{l}f^{t}ac\ot
a^{k}d^{l}f^{t}ba)\delta_{m1}\delta_{n1}\rangle\\
\\
&=&\langle C\ot B, (\sum_{j=0}^{l-1}{a^{k}d^{l-1}f^{t}cr^{-j}s^{t}\ot 
a^{k}d^{l-1}f^{t}br^{-j}s^{-t}})\delta_{m0}\delta_{n0} +
((a^{k+1}d^{l}f^{t} - a^{k}f^{t}bc^{2}(r^{-l}-r^{l}))\\
& &\ot (a^{k+1}d^{l}f^{t}br -
a^{k}f^{t}b^{2}cr(r^{-l}-r^{l}))\delta_{m1}\delta_{n1}\rangle
\end{eqnarray*}

Again, applying the pairings $\langle B, g\rangle$ and $\langle C,
g\rangle$ proves the result $\langle CB, g\rangle$. Now,
\begin{eqnarray*}
r\langle BC, g\rangle - r^{-1}\langle CB, g\rangle =
(r\sum_{j=0}^{k-1}{r^{2(j-l)}} -
r^{-1}\sum_{j=0}^{l-1}{r^{-2j}})\delta_{m0}\delta_{n0}
\end{eqnarray*}

Using the formula 
\begin{eqnarray*}
\sum_{j=0}^{k-1}{x^{j}}= \frac{1-x^{k}}{1-x}
\end{eqnarray*}
we obtain
\[
r\langle BC, g\rangle - r^{-1}\langle CB, g\rangle =
\frac{r^{2(k-l)}-1}{r-r^{-1}}  \delta_{m0}\delta_{n0}
\]
Also obtained as a corollary from the action on
$g=a^{k}d^{l}f^{t}b^{m}c^{n}$ of the monomials in $U_{r,s} $are the
following relations
\[
\begin{array}{lll}
\langle \left[ A, B\right], g\rangle = \langle B, g\rangle,&
\langle \left[ A, C\right], g\rangle =-\langle C, g\rangle&\\
\langle \left[ D, B\right], g\rangle =-\langle B, g\rangle,&
\langle \left[ D, C\right], g\rangle = \langle C, g\rangle&\\
\langle \left[ A, D\right], g\rangle =0,&
\langle \left[ A, F\right], g\rangle =0,&
\langle \left[ D, F\right], g\rangle =0,\\
\langle \left[ B, F\right], g\rangle = \langle B, g\rangle,&
\langle \left[ C, F\right], g\rangle =-\langle C, g\rangle&
\end{array}
\]

From this one obtains the commutation relations in the algebra $U_{r,s}$
dual to $A_{r,s}$ as
\[
rBC - r^{-1}CB = \frac{r^{2(A-D)}-\ident_{U}}{r-r^{-1}}
\]
\[
\begin{array}{lll}
\left[A,B\right]=B, &\left[A,C\right]=-C, &\left[D,B\right]=-B\\
\left[D,C\right]=C, &\left[A,D\right]=0&\\
\left[A,F\right]=0, &\left[D,F\right]=0&\\
\left[B,F\right]=B, &\left[C,F\right]=-C&
\end{array}
\]

The Hopf structure of the dual algebra is obtained and the coproduct of
the elements of the dual is given by
\begin{eqnarray*}
\cop(A) &=& A\ot \ident + \ident\ot A\\
\cop(B) &=& B\ot r^{A-D}s^{-F} + \ident\ot B\\
\cop(C) &=& C\ot r^{A-D}s^{F} + \ident\ot C\\
\cop(D) &=& D\ot \ident + \ident\ot D\\
\cop(F) &=& F\ot \ident + \ident\ot F
\end{eqnarray*}

The counit $\cnt(Y) = 0;$ where $Y = A, B, C, D, F$ and the antipode is
given as
\begin{eqnarray*}
S(A) &=& - A\\
S(B) &=& - Br^{-(A-D)}s^{F}\\
S(C) &=& - Cr^{-(A-D)}s^{-F}\\
S(D) &=& - D\\
S(F) &=& - F
\end{eqnarray*}

\section{The $R$-matrix approach}

Here we obtain the dual algebra employing the $R$-matrix for $\Grs$. We
work with a different $9\times 9$ $R$-matrix by labelling the index $3$ as
$0$ and transposing the $R$-matrix given in [9]. This reads
\[
R=\begin{pmatrix}
r & 0 & 0 & 0\\
0 & S^{-1} & 0 & 0\\
0 & \Lambda & S & 0\\
0 & 0 & 0 & R_{r}
\end{pmatrix}
\]
in block form i.e. in the order $(00)$, $(01)$, $(02)$, $(10)$, $(20)$,
$(11)$, $(12)$, $(21)$, $(22)$ where
\[
R_{r}=\begin{pmatrix}
r & 0 & 0 & 0\\
0 & 1 & 0 & 0\\
0 & \lambda & 1 & 0\\
0 & 0 & 0 & r
\end{pmatrix} \quad
S=\begin{pmatrix}
s & 0\\
0 & 1
\end{pmatrix} ;
\Lambda=\begin{pmatrix}
\lambda & 0\\
0 & \lambda
\end{pmatrix};
\lambda = r-r^{-1}
\]
Also, the the matrix of generators is transposed as
\[
\ct = \begin{pmatrix}
f & 0\\
0 & T
\end{pmatrix} \quad
\text{with} \quad
T=\begin{pmatrix}
a & b\\
c & d
\end{pmatrix}
\]
The linear functionals $(L^{\pm})^a_{b}$ are defined by their value on the
elements of the matrix of generators $T$
\[
\langle (L^{\pm})^a_{b}, T^c_{d} \rangle = (R^{\pm})^{ac}_{bd}
\]
where
\begin{eqnarray*}
(R^{+})^{ac}_{bd} &=& c^{+} (R)^{ca}_{db}\\
(R^{-})^{ac}_{bd} &=& c^{-} (R^{-1})^{ac}_{bd}
\end{eqnarray*}
and $c^{+}$ , $c^{-}$ are free parameters. Matrices $(L^{\pm})^a_{b}$
satisfy
\[
\langle (L^{\pm})^a_{b}, uv \rangle =
\langle (L^{\pm})^a_{c}\ot (L^{\pm})^c_{d}, u\ot v \rangle =
(L^{\pm})^a_{c}(u) (L^{\pm})^c_{d}(v)
\]
\[
\text{i.e.} \quad
\cop (L^{\pm})^a_{b} = (L^{\pm})^a_{c}\ot (L^{\pm})^c_{d}
\]
For $\Grs$, the $(R^{+})$ and $(R^{-})$ matrices read
\[ 
(R^{+})= c^{+} \begin{pmatrix}
r & 0 & 0 & 0\\
0 & S^{-1} & \Lambda & 0\\
0 & 0 & S & 0\\
0 & 0 & 0 & R^{T}_{r}
\end{pmatrix}; \quad
(R^{-})= c^{-} \begin{pmatrix}
r & 0 & 0 & 0\\
0 & S & 0 & 0\\
0 & -\Lambda & S^{-1} & 0\\
0 & 0 & 0 & R^{-1}_{r}
\end{pmatrix}
\]
Note that $R^{-1}_{r}=R_{r^{-1}}$. The $(\cl^{\pm})^a_{b}$ functionals are
dual to the elements $\ct^a_{b}$ in the fundamental representation. From
the above, it is to be noted that the $\cl^{\pm}$ matrices are obtained by
evaluating their action on the generators of the quantum group $\Grs$.
Using the duality pairings, we get

\begin{eqnarray*}
\cl^{+} &=& c^{+}r \begin{pmatrix}
s^{\frac{1}{2} (\tf-H_{2}-1)}r^{\frac{1}{2} (\tf-H_{1}-1)} & 0 & 0\\   
0 & s^{\frac{1}{2}(\tf-H_{1}+1)}r^{\frac{1}{2}(-\tf+H_{2}-1)}
& r^{-1}\lambda \tc\\
0 & 0 & s^{-\frac{1}{2}(\tf+H_{1}-1)}r^{\frac{1}{2}(-\tf-H_{2}-1)}\\
\end{pmatrix}\\
\\
\\
\cl^{-} &=& c^{-}r^{-1} \begin{pmatrix}
s^{-\frac{1}{2} (\tf-H_{2}-1)}r^{-\frac{1}{2} (\tf-H_{1}-1)} & 0 & 0\\
0 & s^{-\frac{1}{2}(\tf-H_{1}+1)}r^{-\frac{1}{2}(-\tf+H_{2}-1)} & 0\\
0 & -r \lambda \tb & 
s^{\frac{1}{2}(\tf+H_{1}-1)}r^{-\frac{1}{2}(-\tf-H_{2}-1)}\\
\end{pmatrix}  
\end{eqnarray*}

where $H_{1}=\ta+\td$, $H_{2}=\ta-\td$ and $\ta, \tb, \tc, \td, \tf$ are
elements of the algebra dual to $\Grs$. More conveniently, one can write

\[
\cl^{+}= \begin{pmatrix}
J & 0 & 0\\
0 & M & P\\
0 & 0 & N\\
\end{pmatrix} \quad \text{and} \quad
\cl^{-}= \begin{pmatrix}
J^{-1} & 0 & 0\\
0 & M^{-1} & 0\\
0 & Q & N^{-1}\\
\end{pmatrix}
\]
where
\[
\begin{array}{ll}
&J=s^{\frac{1}{2} (\tf-H_{2}-1)}r^{\frac{1}{2} (\tf-H_{1}+1)}\\
&M=s^{\frac{1}{2} (\tf-H_{1}+1)}r^{\frac{1}{2} (-\tf+H_{2}+1)}\\
&N=s^{-\frac{1}{2}(\tf+H_{1}-1)}r^{\frac{1}{2} (-\tf-H_{2}+1)}\\
&P=\lambda \tc\\
&Q=-\lambda \tb
\end{array}
\]
These can also be arranged in terms of smaller $L^{+}$ and $L^{-}$
matrices
\[
\begin{array}{llll}
&\cl^{+}=c^{+}\begin{pmatrix}
J & 0\\
0 & L^{+}
\end{pmatrix} \quad &\text{where} \quad
&L^{+}=\begin{pmatrix}
M & P\\
0 & N
\end{pmatrix}\\
&\cl^{-}=c^{-}\begin{pmatrix}
J^{-1} & 0\\
0 & L^{-}
\end{pmatrix} \quad &\text{where} \quad
&L^{-}=\begin{pmatrix}
M^{-1} & 0\\
Q & N^{-1}  
\end{pmatrix}
\end{array}
\]

\subsection*{Commuation relations of the dual}

The algebra dual to $\Grs$ generated by functionals (matrices) $\cl^{\pm}$
satisfy the $q$-commuation relations (the so-called $R\cl\cl$- relations)
\begin{eqnarray*}
R_{12}\cl^{\pm}_{2}\cl^{\pm}_{1} &=& \cl^{\pm}_{1}\cl^{\pm}_{2}R_{12}\\
R_{12}\cl^{+}_{2}\cl^{-}_{1} &=& \cl^{-}_{1}\cl^{+}_{2}R_{12}
\end{eqnarray*}

where $\cl^{\pm}_{1}=\cl^{\pm}\ot \ident$ and  $\cl^{\pm}_{2}=\ident\ot
\cl^{\pm}$. The dual algebra for $\Grs$ can be constructed using the
theory developed in [1]. Since $\Grs$ is a quotient Hopf
algebra, the $R$-matrix for the $R\cl\cl$ - relations is different from
the one used in the $RTT$ - relations for the $\Grs$ algebra. Instead, it
is necessary to amend the $R$-matrix to remove relations which are
inconsistent with the quotient structure. For $\Grs$, this means that the
$R\cl\cl$ - relations are constructed with the $R$-matrix

\[
R_{12}=c^{-}{\langle \cl^{-},\ct \rangle}^{-1}
=\begin{pmatrix}
r & & &\\
& S^{-1} & &\\
& & S &\\
& & & R_{r}
\end{pmatrix}
\]
Evaluating $\cl^{\pm}_{1}$, $\cl^{\pm}_{2}$ matrices and substituting in
the above $R\cl\cl$- relations yields the dual algebra commutation 
relations.
From
$R_{12}\cl^{-}_{2}\cl^{-}_{1} = \cl^{-}_{1}\cl^{-}_{2}R_{12}$ and
$R_{12}\cl^{+}_{2}\cl^{+}_{1} = \cl^{+}_{1}\cl^{+}_{2}R_{12}$
we obtain
\begin{eqnarray*}
R_{r}L^{-}_{2}L^{-}_{1} &=& L^{-}_{1}L^{-}_{2}R_{r}\\
R_{r}L^{+}_{2}L^{+}_{1} &=& L^{+}_{1}L^{+}_{2}R_{r}\\
MJ &=& JM\\
NJ &=& JN\\
JQ &=& s^{-1}QJ\\
JP &=& sPJ
\end{eqnarray*}
where
\[
\begin{array}{llllll}
&R_{r}L^{-}_{2}L^{-}_{1}=L^{-}_{1}L^{-}_{2}R_{r} &\Longrightarrow
&MQ=rQM,\quad QN=rNQ &\text{and} &NM=MN\\
&R_{r}L^{+}_{2}L^{+}_{1}=L^{+}_{1}L^{+}_{2}R_{r} &\Longrightarrow
&PM=rMP,\quad NP=rPN &\text{and} &NM=MN
\end{array}
\]
In addition, the cross relation 
$R_{12}\cl^{+}_{2}\cl^{-}_{1} = \cl^{-}_{1}\cl^{+}_{2}R_{12}$ yields
$R_{r}L^{+}_{2}L^{-}_{1} = L^{-}_{1}L^{+}_{2}R_{r}$ which further
implies
\[
QP-PQ=-\lambda(N^{-1}M-NM^{-1})
\]
Simplifying the above, we get the following commutation relations
\[
\begin{array}{lll}
&[\ta,\tb]=\tb, \quad &[\ta,\tc]=-\tc\\
&[\td,\tb]=-\tb, \quad &[\td,\tc]=\tc\\
&[\ta,\td]=0, \quad &[\tf, \bullet]=0\\
\end{array}
\]
and 
\[
[\tb,\tc]=
\frac{(r^{\ta-\td}s^{\tf})-{(r^{\ta-\td}s^{\tf})}^{-1}}{r-r^{-1}}= 
\frac{r^{\ta-\td+\gamma \tf}-r^{-(\ta-\td+\gamma \tf)}}{r-r^{-1}} \quad
\text{where} \quad
\gamma =\frac{\ln s}{\ln r}
\]
So, defining $H=\ta-\td+\gamma \tf$, $X_{+}=\tb$ and $X_{-}=\tc$ we obtain
\[
[H,X_{\pm}]=2X_{\pm}; \quad [X_{+},X_{-}]=[H]; \quad [\tf,\bullet]=0
\]
i.e. the Drinfeld-Jimbo form of a single parameter deformation of
$gl(2)\oplus gl(1)$.

\section{Differential Calculus on $\Grs$}

We now proceed towards the construction of a differential calculus on the
$\Grs$ quantum group. We use Jurco's constructive procedure [10] based on
the $R$-matrix formulation and using the $\ct$- and the $\cl^{\pm}$
matrices of the previous section.

\subsection*{Quantum one-forms}

Define $\omega$ to be the basis of all left-invariant quantum one-forms.
So, we have
\[
\cop_{L}(\omega)= \ident\ot \omega
\]
This defines the left action on the bimodule $\Gamma$ (space of quantum
one-forms). The bimodule $\Gamma$ is further characterised by the
commutation relations between $\omega$ and $a\in A$ ($A$ being the $\Grs$
Hopf algebra),
\[
\omega a = (f\ast a) \omega
\]
The left convolution product is defined
\[
f\ast a = (\ident\ot f) \cop (a)
\]
where $f\in A'(=Hom(A, \cal C))$ i.e. belongs to the dual. This means
\[
\omega a = (\ident\ot f) \cop (a) \omega
\]
Now, the linear functional $f$ is defined in terms of the $\cl^{\pm}$
matrices as
\[
f = S(\cl^{+})\cl^{-}
\]
Thus, we have
\[
\omega a = [(\ident\ot S(\cl^{+})\cl^{-}) \cop (a)] \omega
\]
In terms of components, one can write
\[
\omega_{ij} a = [(\ident\ot S(l^{+}_{ki})l^{-}_{jl}) \cop (a)] \omega_{kl}
\]
using the expressions $\cl^{\pm} = l^{\pm}_{ij}$ and $\omega =
\omega_{ij}$ where $i, j = 1 .. 3$.\par

For $\Gamma$ to be a bicovariant bimodule, the right coaction is given as
\[
\cop_{R} (\omega) = \omega\ot M
\]
where functionals $M$ are defined in terms of the matrix of generators
$T$,
\[
M = TS(T)
\]
Again, in component form, one can write
\[
\cop_{R} (\omega_{ij}) = \omega_{kl}\ot t_{ki}S(t_{jl})
\]
Using the above formulae, we obtain the commutation relations of all the
left-invariant one forms with the elements of the $\Grs$ quantum group as
follows
\[
\begin{alignat*}{2}
\omega^{0} a &= r^{2}s^{2}a\omega^{0} & \qquad \qquad
\omega^{0} b &= r^{2}b\omega^{0}\\
\omega^{1} a &= r^{-2}a\omega^{1} & \qquad \qquad
\omega^{1} b &= b\omega^{1}\\
\omega^{+} a &= r^{-1}a\omega^{+} & \qquad \qquad
\omega^{+} b &= r^{-1}b\omega^{+} - \lambda r^{-1}a\omega^{1}\\
\omega^{-} a &= r^{-1}a\omega^{-} -\lambda r^{-1}b\omega^{1} & \qquad
\qquad
\omega^{-} b &= r^{-1}b\omega^{-}\\
\omega^{2} a &= a\omega^{2} -\lambda b\omega^{+} & \qquad \qquad
\omega^{2} b &= r^{-2}b\omega^{2} -\lambda r^{-1}a\omega^{-} +
\lambda^{2} b\omega^{1}\\
\omega^{0} c &= r^{2}s^{2}c\omega^{0} & \qquad \qquad
\omega^{0} d &= r^{2}d\omega^{0}\\
\omega^{1} c &= r^{-2}c\omega^{1} & \qquad \qquad
\omega^{1} d &= d\omega^{1}\\
\omega^{+} c &= r^{-1}c\omega^{+} & \qquad \qquad
\omega^{+} d &= r^{-1}d\omega^{+} - \lambda r^{-1}c\omega^{1}\\
\omega^{-} c &= r^{-1}c\omega^{-} -\lambda r^{-1}d\omega^{1} & \qquad
\qquad
\omega^{-} d &= r^{-1}d\omega^{-}\\
\omega^{2} c &= c\omega^{2} -\lambda d\omega^{+} & \qquad \qquad
\omega^{2} d &= r^{-2}d\omega^{2} -\lambda r^{-1}c\omega^{-} +
\lambda^{2} d\omega^{1}
\end{alignat*}
\]
\begin{align*}
\omega^{0} f &= f\omega^{0}\\
\omega^{1} f &= s^{-2}f\omega^{1}\\
\omega^{+} f &= s^{-1}f\omega^{+}\\
\omega^{-} f &= s^{-1}f\omega^{-}\\
\omega^{2} f &= f\omega^{2}
\end{align*}

where $\omega^{0}=\omega_{11}, \omega^{1}=\omega_{22},
\omega^{+}=\omega_{23}, \omega^{-}=\omega_{32}, \omega^{2}=\omega_{33}$
and the components $\omega_{12}, \omega_{13}, \omega_{21}, \omega_{31}$
have null contribution given the structure of the $\ct$- matrix (i.e.
$t_{12}=t_{13}=t_{21}=t_{31}=0$).

\subsection*{Vector Fields}

The linear space $\Gamma$ (space of all left invariant
one-forms) contains a bi-invariant element $\tau = \sum_{i}\omega_{ii}$
which can be used to define a derivative on $A$ (the $\Grs$ Hopf algebra).
For $a\in A$, one sets
\[
\mathbf{d} a = \tau a - a \tau
\]
\[ \text{Now} \qquad
\omega_{ii} a = [(\ident\ot S(l^{+}_{ki})l^{-}_{il}) \cop (a)] \omega_{kl}
\]
\[ \text{So} \qquad
\mathbf{d}a = [(\ident\ot \chi_{kl}) \cop (a)] \omega_{kl}
\]
where $\chi_{kl} = S(l^{+}_{ki})l^{-}_{il} - \delta_{kl}\cnt$, $\cnt$
being the counit. Denote
\[
\chi_{ij} = S(l^{+}_{ik})l^{-}_{kj} - \delta_{ij}\cnt
\]
or more compactly
\[
\chi = S(\cl^{+})\cl^{-} - \ident\cnt
\]
the matrix of left-invariant vector fields $\chi_{ij}$ on $A$. On elements
of $\Grs$ the vector fields act as
\begin{eqnarray*}
\chi_{ij} a &=& (S(l^{+}_{ik})l^{-}_{kj} - \delta_{ij}\cnt)a\\
\chi_{ij} a &=& \langle S(l^{+}_{ik})l^{-}_{kj}, a\rangle -
\delta_{ij}\cnt (a); \qquad \qquad a\in \Grs
\end{eqnarray*}
We obtain explicitly the following

\begin{alignat*}{2}
\chi_{0}(a) &= r^{2}s^{2}-1 & \qquad \chi_{0}(b) &= 0\\
\chi_{1}(a) &= r^{-2}-1 & \qquad \chi_{1}(b) &= 0\\
\chi_{+}(a) &= 0 & \qquad \chi_{+}(b) &= 0\\
\chi_{-}(a) &= 0 & \qquad \chi_{-}(b) &= -(r-r^{-1})\\
\chi_{2}(a) &= 0 & \qquad \chi_{2}(b) &= 0\\
\\
\chi_{0}(c) &= 0 & \qquad \chi_{0}(d) &= r^{2}-1\\
\chi_{1}(c) &= 0 & \qquad \chi_{1}(d) &= (r-r^{-1})^{2}\\
\chi_{+}(c) &= -(r-r^{-1}) & \qquad \chi_{+}(d) &= 0\\
\chi_{-}(c) &= 0 & \qquad \chi_{-}(d) &= 0\\
\chi_{2}(c) &= 0 & \qquad \chi_{2}(d) &= r^{-2}-1
\end{alignat*}
\begin{align*}
\chi_{0}(f) &= 0\\
\chi_{1}(f) &= s^{-2}-1\\
\chi_{+}(f) &= 0\\
\chi_{-}(f) &= 0\\
\chi_{2}(f) &= 0
\end{align*}

where $\chi_{0}=\chi_{11}, \chi_{1}=\chi_{22}, \chi_{+}=\chi_{23},
\chi_{-}=\chi_{32}, \chi_{2}=\chi_{33}$ and again (by previous argument)
the components $\chi_{12}, \chi_{13}, \chi_{21}, \chi_{31}$ have null
contribution. The left convolution products are given as

\begin{alignat*}{2}
\chi_{0} \ast a &= ((rs)^{2}-1)a & \qquad \chi_{0} \ast b &= (r^{2}-1)b\\
\chi_{1} \ast a &= (r^{-2}-1)a & \qquad \chi_{1} \ast b &=
((r-r^{-1})^{2})b\\
\chi_{+} \ast a &= -(r-r^{-1})b & \qquad \chi_{+} \ast b &= 0\\
\chi_{-} \ast a &= 0 & \qquad \chi_{-} \ast b &= -(r-r^{-1})a\\
\chi_{2} \ast a &= 0 & \qquad \chi_{2} \ast b &= (r^{-2}-1)b\\
\\
\chi_{0} \ast c &= ((rs)^{2}-1)c & \qquad \chi_{0} \ast d &= (r^{2}-1)d\\
\chi_{1} \ast c &= (r^{-2}-1)c & \qquad \chi_{1} \ast d &=
((r-r^{-1})^{2})d\\
\chi_{+} \ast c &= -(r-r^{-1})d & \qquad \chi_{+} \ast d &= 0\\
\chi_{-} \ast c &= 0 & \qquad \chi_{-} \ast d &= -(r-r^{-1})c\\
\chi_{2} \ast c &= 0 & \qquad \chi_{2} \ast d &= (r^{-2}-1)d
\end{alignat*}
\begin{align*}   
\chi_{0} \ast f &= 0\\
\chi_{1} \ast f &= (s^{-2}-1)f\\
\chi_{+} \ast f &= 0\\
\chi_{-} \ast f &= 0\\
\chi_{2} \ast f &= 0
\end{align*}

\subsection*{Exterior Derivatives}

Using the formula $\mathbf{d} a =\sum_{i}(\chi_{i} \ast a) \omega^{i}$ for
$a \in A$, we obtain the action of the exterior derivative on the
generating elements of $\Grs$

\begin{eqnarray*}
\mathbf{d} a &=& ((rs)^{2}-1)a\omega^{0} + (r^{-2}-1)a\omega^{1} - \lambda
b\omega^{+}\\
\mathbf{d} b &=& (r^{2}-1)b\omega^{0} + \lambda^{2}b\omega^{1} - \lambda
a\omega_{-} + (r^{-2}-1)b\omega^{2}\\
\mathbf{d} c &=& ((rs)^{2}-1)c\omega^{0} + (r^{-2}-1)c\omega^{1} - \lambda
d\omega^{+}\\
\mathbf{d} d &=& (r^{2}-1)d\omega^{0} + \lambda^{2}d\omega^{1} - \lambda  
c\omega_{-} + (r^{-2}-1)d\omega^{2}\\
\mathbf{d} f &=& (s^{-2}-1)f\omega^{1}
\end{eqnarray*}

where $\lambda=r-r^{-1}$. The exterior derivative $\mathbf{d}: A \lra
\Gamma$ satisfies the  Leibnitz rule and $\mathbf{d}A$ generates $\Gamma$
as a left $A$-module. This then defines a first order differential
calclulus $(\Gamma, \mathbf{d})$ on $\Grs$. Furthermore, the calculus is
bicovariant due to the coexistence of the left and the right actions
\begin{eqnarray*}
\cop_{L}:\Gamma \lra A \ot \Gamma\\
\cop_{R}:\Gamma \lra \Gamma \ot A
\end{eqnarray*}
since $\mathbf{d}$ has the invariance property
\begin{eqnarray*}
\cop_{L}\mathbf{d}=(\ident \ot \mathbf{d})\cop\\
\cop_{R}\mathbf{d}=(\mathbf{d} \ot \ident)\cop       
\end{eqnarray*}

Motivated by the observation that $\Grs$ provides a realisation of the two
parameter quantum group $\GLpqtwo$, we wish to investigate the
relationship between the differential calculus on $\Grs$ and that of
$\GLpqtwo$. Again, if we denote the primed generators as those belonging
to $\GLpqtwo$ while the unprimed ones as the generators of $\Grs$, then
the
relation
\[
\mathbf{d} \begin{pmatrix}a'&b'\\c'&d'\end{pmatrix}
=\mathbf{d} \left( f
\begin{pmatrix}a&b\\c&d\end{pmatrix} \right)
\]
provides a relisation of the differential calculus on $\GLpqtwo$ [13],
with the
defining relations between the two sets of deformation parameters
$(p,q)$ and $(r,s)$ as before. In establishing the above homomorphism at
the differential calculus level, one has to ignore the term involving
$\omega^{0}$ as this is the one-form related to the fifth element $f$,
which is nonexistent in the case of $\GLpqtwo$. A large part of the
differential calculus depends on the deformation parameter $r$ with the
second parameter $s$ being related to the element $f$ and corrsponding
one-forms and vector fields.\par

We note that the calculus on $\Grs$ contains the calculus on $\GLqtwo$ and
our results match with those given in [11]. We also expect that the
calculus on $\Grs$ could be obtained by projection from the calculus on
multiparameter $q$-deformed $GL(3)$ [12]. The physical interest in
studing $\Grs$ lies in the observation that when endowed with a $\ast$-
structure, this quantum group would specialise to a two parameter quantum
deformation of $SU(2) \ot U(1)$ which is precisely the gauge group for the
theory of electroweak interactions. Since gauge theories have an obvious
differential geometric description, the above study of differential
calculus could provide insights in constructing a $q$-gauge theory based
on $\Grs$. The homomorphism between $\Grs$ and $\GLpqtwo$ at the
level of differential calculus is then, also expected to play a
significant role.

\section{Conclusions}

In this work, we have explicitly constructed the dual algebra to the
quantum group $\Grs$ and exhibited its Hopf algebraic structure. This has
been obtained using the method described by Sudbery [6] as well as the
$R$-matrix approach [1]. It is to be noted that the algebra sector of the
dual to $\Grs$ is a single parameter deformation of $gl(2)\oplus gl(1)$
whereas the second parameter only appears in the coalgebra. This can be
clearly understood as an application of the twisting procedure. As a next
stage of this analysis, we have constructed a bicovariant differential
calculus on the quantum group $\Grs$, which would be of significance in
constructing a corresponding gauge theory. Furthermore, our analysis shows
that the quantum group homomorphism between $\Grs$ and  $\GLpqtwo$ has a
natural extention at the level of their differential calculii.
\par

\section*{Acknowledgments}

The authors would like to thank Professors Allan Solomon and Tony Sudbery
for enlightening discussions. One of us (D.P.) is also grateful to
Professor Tetsuji Miwa for the kind hospitality at the Research Institute
for Mathematical Sciences (RIMS), Kyoto via a JSPS fellowship, where this
work was initiated.
\par

\section*{References}

[1]	L. D. Faddeev, N. Yu. Reshetikhin and L. A. Takhtajan, {\sl Len.
Math. J.}{\bf 1}, 193 (1990).\par
[2]	Yu I. Manin, Montreal University Preprint CRM - 1561, 1988.\par
[3]	Yu I. Manin, {\sl Comm. Math. Phys.} {\bf 123}, 163 (1989).\par
[4]	S. L. Woronowicz, {\sl Comm. Math. Phys.} {\bf 111}, 613
(1987).\par
[5]	S. L. Woronowicz, {\sl Comm. Math. Phys.} {\bf 122}, 125
(1989).\par
[6]	A. Sudbery, Proc. Workshop on Quantum Groups, Argonne (1990)
eds. T. Curtright, D. Fairlie 	and C. Zachos, pp. 33-51.\par
[7]	V. K. Dobrev, {\sl J. Math. Phys.} {\bf 33}, 3419 (1992).\par
[8]	V. K. Dobrev and P. Parashar, {\sl J. Phys.} {\bf A26}, 6991
(1993).\par
[9]	B. Basu-Mallick, hep-th/9402142; {\sl Intl. J. Mod. Phys.} {\bf
A10}, 2851 (1995).\par
[10]	B. Jurco, {\sl Lett. Math. Phys.} {\bf 22}, 177 (1991); Preprint
CERN-TH 9417/94 (1994).\par
[11]	P. Aschieri and L. Castellani, {\sl Intl. J. Mod. Phys.} {\bf A8},
1667 (1993).\par
[12]	P. Aschieri and L. Castellani, {\sl Phys. Lett.} {\bf B293}, 299
(1992).\par
[13]	F. Muller-Hoissen, {\sl J. Phys.} {\bf A25}, 1703 (1992); X-D Sun
and S-K Wang, {\sl J. Math. Phys.} {\bf 33}, 3313 (1992).\par

\end{document}